\documentclass[11pt,a4paper]{article}
\usepackage[latin1]{inputenc}
\usepackage{a4wide}
\usepackage{amsmath}
\usepackage{amsthm}
\usepackage{amssymb}  
\usepackage{amsfonts}
\usepackage{cite}
\newtheorem{theorem}{Theorem}
\newtheorem{lemma}{Lemma}
\newtheorem{corollary}{Corollary}

\title{Almost complex structures on real Lie supergroups}
\author{Matthias Kalus\footnote{Research supported by the SFB/TR 12, Symmetry and Universality in Mesoscopic Systems, of the Deutsche Forschungsgemeinschaft.}}
\date{}
\begin{document}
\maketitle
%--------------------------------------------------------------------------
% ABSTRACT
%--------------------------------------------------------------------------
\begin{abstract}
\noindent
A complex Lie supergroup can be described as a real Lie supergroup with integrable almost complex structure.  The necessary and sufficient conditions on an almost complex structure on a real Lie supergroup for defining a complex Lie supergroup are deduced. The classification of real Lie supergroups with such almost complex structures yields a new  approach to the known classification of complex Lie supergroups by complex Harish-Chandra superpairs. A universal complexification of a real Lie supergroup is constructed.

\end{abstract}
\setcounter{tocdepth}{1}
%\tableofcontents \vspace*{0.5cm}
%--------------------------------------------------------------------------
% INTRODUCTION
%--------------------------------------------------------------------------

\noindent
The local differential operators on a real Lie supergroup have the structure of a Lie-Hopf superalgebra which can be algebraically constructed from a real Harish-Chandra superpair (see \cite{K}). Conversely starting from a real Harish-Chandra superpair, Kostant constructed in \cite{K} a sheaf of superfunctions by dualising the associated Lie-Hopf superalgebra. This yields a real Lie supergroup and hence an equivalence of categories from real Harish-Chandra superpairs to real Lie supergroups. 

\bigskip\noindent
A construction of complex Lie supergroups from complex Harish-Chandra superpairs using analytic continuation on Grassmann variables was given by Berezin (see \cite{Ber}). Vishnyakova gave a rigorous proof of the equivalence of categories of complex Harish-Chandra superpairs and complex Lie supergroups (see \cite{Vish}). 

\bigskip\noindent
In this article complex Lie supergroups and complex Harish-Chandra superpairs are analyzed as real objects with integrable almost complex structure. This  approach agrees with the known equivalence of categories, but it should be remarked in the literature. As an application of the equivalence of categories, the universal complexification of a real Lie supergroup is constructed. More Details can be found in \cite{Dis}.

\bigskip\noindent
\textit{Contents.} We associate a real supermanifold with integrable almost complex structure $J$ to a complex supermanifold, in particular to a  Lie supergroup (see also e.g. \cite{DM}). The graded version of the Newlander-Nierenberg theorem (see \cite{McH},\cite{Va}) yields the way back on the level of supermanifolds. For forming a complex  Lie supergroup starting from a real one with almost complex structure, $J$ has to satisfy more conditions which we explicitly deduce. The correspondence of complex Lie supergroups and Harish-Chandra superpairs then follows from the real case in \cite{K}. Furthermore the sheaf of holomorphic superfunctions can be expressed by a dualization of the complex  Lie-Hopf superalgebra analogously to the real case. Existence of a universal complexification of a real Lie supergroup with underlying real analytic Lie group is finally derived from the non-graded case (see \cite{Hoch}) in the language of complex Harish-Chandra superpairs.

\bigskip\noindent
\textit{Acknowledgment.} The author wishes to thank A. Huckleberry for his advise and discussions.

\subsubsection* {From complex to real}

Let $\mathcal U=(U, \mathcal O_U \otimes_\mathbb C \Lambda (\mathbb C^n)^\ast)$ for open $U \subset \mathbb C^m\cong \mathbb R^{2m}$ be a complex superdomain as it is used for local coordinate charts of supermanifolds. A holomorphic superfunction on $\mathcal U$ can be regarded as a finite sum of $\mathbb C$-multilinear alternating forms on $\mathbb C^n$ with values in $\mathcal O_U$. Note that an element in $\Lambda^k (\mathbb C^n)^\ast$ is naturally also an $\mathbb R$-multilinear $\mathbb C$-valued $k$-form on $\mathbb R^{2n}$ so an element in $\Lambda^k (\mathbb R^{2n})^\ast \otimes_\mathbb R \mathbb C$. We denote the sheaf of smooth functions with values in $\mathbb K=\mathbb R$, resp. $\mathbb C$ by $\mathcal C^\infty_{\mathbb K,U}$. Embedding holomorphic functions $\mathcal O_U \hookrightarrow \mathcal C^{\infty}_{\mathbb C,U}= \mathcal C^\infty_{\mathbb R,U} \otimes_\mathbb R \mathbb C$ we have 
$$\mathcal A_\mathcal U:=\mathcal O_U \otimes_\mathbb C \Lambda (\mathbb C^{n})^\ast\ \ \subset \ \ \mathcal C_\mathcal U^\mathbb C:=\mathcal C^{\infty}_{\mathbb C,U} \otimes_\mathbb R  \Lambda (\mathbb R^{2n})^\ast \ .$$ 
Furthermore we obtain a real superdomain $\hat{\mathcal U}=(U,\mathcal C_\mathcal U^\mathbb R)$ via the subsheaf $\mathcal C_\mathcal U^\mathbb R:=\mathcal C^\infty_{\mathbb R,U} \otimes_\mathbb R \Lambda (\mathbb R^{2n})^\ast$ of $\mathcal C^\mathbb C_\mathcal U$. We identify respectively the sheaves of superderivations on $\mathcal A_\mathcal U$ with $\mathcal A_\mathcal U \otimes_\mathbb C \mathbb C^{m+n}$ and on $\mathcal C^\mathbb K_\mathcal U$ with $\mathcal C^\mathbb K_\mathcal U \otimes_\mathbb R \mathbb R^{2(m+n)}$  in the canonical way. On $Der(\mathcal C_U^\mathbb R)$ we have the almost complex structure $J_\mathcal U$ (i.e. a $\mathcal C_\mathcal U^\mathbb R$-linear (even) automorphism of $Der(\mathcal C_\mathcal U^\mathbb R)$ satisfying $J_\mathcal U^2=-Id$) given by multiplication with $i$ on the identified $\mathbb R^{2(m+n)}=\mathbb C^{(m+n)}$. The elements in the $-i$ eigenspace of its $\mathbb C$-linear continuation to $Der(\mathcal C^\mathbb C_\mathcal U)$ have the common kernel $\mathcal A_\mathcal U$. We denote the common kernel of the $+i$ eigenvectors by $\overline{\mathcal A}_\mathcal U$ and have the real isomorphism of sheaves $\psi_\mathcal U:\mathcal A_\mathcal U \to \overline{\mathcal A}_\mathcal U$, $f \mapsto \overline{f}$ defined by complex conjugation.

\bigskip\noindent
Let now $\Phi=(\varphi,\Phi^\ast):\mathcal U \to \mathcal V$ be a morphism of complex superdomains as it appears as the  gluing of coordinate charts in an atlas of a complex supermanifold or as the local version of a morphism between complex supermanifolds. Its pullback $\Phi^\ast:\mathcal A_\mathcal V \to \Phi_\ast\mathcal A_\mathcal U$ and the associated morphism $\psi_\mathcal U\circ \Phi^\ast \circ \psi_\mathcal V:\overline{\mathcal A}_\mathcal V \to \Phi_\ast\overline{\mathcal A}_\mathcal U$ can be completed to a morphism on $\mathcal C^\mathbb C_\mathcal V $ and restricted to a map $\hat\Phi^\ast:\mathcal C_\mathcal V^\mathbb R \to \hat\Phi_\ast \mathcal C_\mathcal U^\mathbb R$. Hence we naturally obtain a morphism of real superdomains $\hat\Phi=(\hat\varphi,\hat{\Phi}^\ast):\hat{\mathcal U} \to \hat{\mathcal V}$. It is compatible with the almost complex structures in the sense of $\hat\Phi_\ast \circ J_\mathcal U=J_\mathcal V \circ \hat\Phi_\ast$. So applied to the gluing of coordinate charts, this construction associates a real supermanifold with integrable almost complex structure to a complex supermanifold. On global level it associates  a morphism of real supermanifolds compatible with the almost complex structure to a morphism of complex supermanifolds. Hence we obtain from a complex Lie supergroup by this construction a real Lie supergroup with integrable almost complex structure and compatible multiplication, inverse and unity morphisms. A morphism of complex Lie supergroups yields a morphism of real Lie supergroups compatible with the almost complex structures.

\subsubsection* {From real to complex}

Real Lie supergroups can be identified with Harish-Chandra superpairs, i.e. pairs $(G,\mathfrak g)$ consisting of a real Lie group $G$ and a real Lie superalgebra $\mathfrak g=\mathfrak g_{\bar 0} \oplus\mathfrak g_{\bar 1}$ such that $\mathfrak g_{\bar 0}=Lie(G)$ and the representation of $\mathfrak g_{\bar 0}$ on $\mathfrak g_{\bar 1}$ integrates to a representation of $G$ (see \cite{K}). Furthermore Kostant shows in \cite{K} that the supermanifold structure on a real Lie supergroup $\mathcal G=(G,\mathcal C_\mathcal G^\mathbb R)$ and the morphisms of multiplication, inverse and unity induce the structure of a Lie-Hopf superalgebra on 
\begin{align*}
  {\mathbf  C}_\mathcal G^\ast:= 
\big\{\varphi \in &Hom_{\mathbb R-vect}({\mathcal C_\mathcal G^\mathbb R}(G),\mathbb R)\big|  \\ & \quad
\exists \mbox{ ideal }  I \subset {\mathcal C_\mathcal G^\mathbb R}(G) \mbox{ with } 
codim(I)<\infty,I \subset ker(\varphi) \big\} \ 
\end{align*}
It is shown in \cite {K} that ${\mathbf  C}_\mathcal G^\ast$ can be identified with $\mathbb R(G)\# E(\mathfrak g)$. The first factor is the group ring of $G$, the  second is the universal enveloping algebra of $\mathfrak g$ and $\#$ denotes the semidirect product on the tensor product of vector spaces by
$$
(x_1\# y_1)\cdot (x_2\# y_2)=\left(x_1 \cdot x_2\right) \# (  {Ad}(x_2^{-1})(y_1) \cdot y_2) \quad \mbox{for }  x_1,x_2 \in G,\ y_1,y_2\in E(\mathfrak g) \ .
$$
Here $g\#X$ is interpreted as applying the operator $X$ and afterward evaluating  at $g$.  On the other hand starting with a Harish-Chandra superpair and the associated Lie-Hopf superalgebra $\mathbb R(G)\# E(\mathfrak g)$, Kostant constructs the sheaf of superfunctions on the associated real Lie supergroup by 
\begin{align}\label{eqeq}
\begin{array}{rl}{\mathcal C}_\mathcal G^\mathbb K(U)\ \ = &
\Big\{\Phi \in Hom_{\mathbb R-vect} 
 (\mathbb R(U) \# E(\mathfrak g), \mathbb K) \ \big| \\ 
& \qquad\Big(U \to \mathbb K, \ g \mapsto \Phi(g \# Z)\Big)  
\in \mathcal C^\infty_{\mathbb K, G}(U) \  \forall  Z \in E(\mathfrak g) \Big\}  \quad \mbox{for } \mathbb K=\mathbb R\ .\end{array}
\end{align}  
This direct correspondence of Harish-Chandra superpairs and Lie supergroups is highly dependent on the softness of the sheaf of smooth functions and can not be transported directly to the complex setting. However  it is possible, as will be shown in the following.

\bigskip\noindent
Let $J$ be an almost complex structure on a real supermanifold $\mathcal M$ (later on a real Lie supergroup $\mathcal G$). The necessary and sufficient condition on $J$ to be integrable, i.e. to define a complex supermanifold structure, is that the supercommutator on vector-fields continued $\mathbb C$-linearly to $Der(\mathcal C^\mathbb C_\mathcal M)$ preserves the eigenspaces of $J$. (See \cite{McH},\cite{Va} for the general result and parallel arguments to \cite{KN} chap.~IX.2 for equivalence to the condition stated here.) A morphism of real supermanifolds with integrable almost complex structures becomes a morphism of complex supermanifolds if and only if it is compatible with the almost complex structures.  We aim at further conditions on an integrable $J$ on a real Lie supergroup $\mathcal G$ to define a complex Lie supergroup, i.e.~the conditions coming from $J$-compatibility of multiplication, inverse and unity morphisms.

\bigskip\noindent
Let now $\mathcal G=(G,\mathcal C_\mathcal G^\mathbb R)$ be a real Lie supergroup associated to the superpair $(G,\mathfrak g)$ and let $J$ be an integrable almost complex structure on the supermanifold $\mathcal G$. For $X \in \mathfrak g$ regarded as a left-invariant derivation on ${\mathcal C}_\mathcal G^\mathbb R$, we obtain for any $g\in G$ a well-defined element $J_g(X) \in \mathfrak g$ such that $g\# J_g(X)=J(g\#X)$ as derivations. The condition that the multiplication map $(m,m^\ast)$ on $\mathcal G$ is a morphism of complex supermanifolds, i.e. preserves $J$, is translated to the Lie-Hopf superalgebra as $ (h\#1)\cdot J(g\#X)=J((h\#1)\cdot(g\#X))$ for all $g,h \in G$ and $X \in \mathfrak g$ which yields with $h=g^{-1}$ that $J_g(X)=J_e(X)$.  So $J$ is supposed to map left-invariant derivations to left-invariant derivations and hence to restrict to a map $J_\mathfrak g:\mathfrak g \to \mathfrak g$. Furthermore for homogeneous $X_\pm$ in the $\pm i$ eigenspace of $J$, use on the left hand side of
$$(e\#X_+)\cdot (e\#X_-)-(-1)^{|X_+||X_-|}(e\#X_-)\cdot (e\#X_+)=e\#[X_+,X_-]$$
the compatibility with $J$ first in the first, then in the second argument of $(m,m^\ast)$. This yields identical results with different signs. So $[X_+,X_-]$ vanishes. Together with the graded Newlander-Nierenberg theorem we obtain $J$-linearity in both arguments of the superbracket. These conditions are already sufficient:
\begin{theorem}
A real Lie supergroup $\mathcal G$ with almost complex structure $J$ induces a complex Lie supergroup if and only if $J$  preserves left-invariance of superderivations and the Lie superbracket is $J$-linear in both arguments, i.e.~$J$ comes from a complex structure on the Lie superalgebra $\mathfrak g$.
\end{theorem}
\begin{proof} If $J$ satisfies the conditions above, then it is integrable due to the graded version of the Newlander-Nierenberg theorem. Furthermore $J$ is compatible with the adjoint action of $G$ on $\mathfrak g$, so $J$ can be continued to $\mathbb R(G) \# E(\mathfrak g)$ compatible with multiplication, inverse $(g\#X \mapsto -g^{-1}\#Ad(g)(X)$ for $g \in G$ and $X\in \mathfrak g$) and unity. This includes that the corresponding morphisms are morphisms of complex supermanifolds. 
\end{proof}

\noindent
From the first section and the fact, that a complex Harish-Chandra superpair is a real Harish-Chandra superpair with a complex structure on the Lie superalgebra $\mathfrak g$, we can follow a result first proved in \cite{Vish}:
\begin{corollary}
  The category of complex Lie supergroups is equivalent to the category of complex Harish-Chandra superpairs.
\end{corollary}

\noindent
The sheaf of holomorphic superfunctions for a complex Lie supergroup can be constructed from a Harish-Chandra superpair parallel to (\ref{eqeq}): Let $(G,\mathfrak g)$ be a real Harish-Chandra superpair with almost complex structure $J$ as in the theorem. Note that the $+i$ eigenspace of $J$ in $\mathfrak g \otimes_\mathbb R \mathbb C$ consists of the elements  $X-iJ(X)$, $X\in \mathfrak g$ and hence is canonically isomorphic to $\mathfrak g$. Furthermore $\Phi \in \mathcal C^\mathbb C_\mathcal G(U)$ is holomorphic if and only if it vanishes on all monomials in $\mathbb R(U)\# E(\mathfrak g)$ containing an element in the $-i$ eigenspace of $J$. We obtain:
\begin{corollary}
The sheaf of holomorphic superfunctions on a complex Lie supergroup $\mathcal G$ associated to the complex superpair $(G,\mathfrak g)$ can be explicitly given by:
\begin{align*}
{\mathcal A}_\mathcal G(U)= 
\Big\{ \Phi \in Hom_{\mathbb C-vect} 
& (\mathbb C(U) \# E(\mathfrak g ), \mathbb C) \ \big| \\ 
& \Big(U \to \mathbb C, \ g \mapsto \Phi(g \# Z)\Big)  
\in \mathcal O_{G}(U) \  \forall  Z \in E(\mathfrak g) \Big\} \ 
\end{align*}
on open subsets $U\subset G$.
\end{corollary}

\subsubsection* {Universal complexification}
\noindent
Let now $G$ be a real analytic group. We finally prove existence of a universal complexification of a real Lie supergroup with underlying Lie group $G$. A universal complexification of a real Lie supergroup $\mathcal G$ is a complex Lie supergroup $\mathcal G^\mathbb C$ and a morphism of real Lie supergroups $\Gamma:\mathcal G \to \mathcal G^\mathbb C$ with the universal property: for any morphism of real Lie supergroups $\Phi:\mathcal G\to \mathcal H$ into a complex Lie supergroup $\mathcal H$ there exists a unique morphism of complex Lie supergroups $\Phi^\mathbb C:\mathcal G^\mathbb C \to \mathcal H$ such that $\Phi^\mathbb C\circ \Gamma=\Phi$. Note that existence of a universal complexification includes uniqueness up to isomorphisms of complex Lie supergroups. 

\bigskip\noindent
For real analytic Lie groups $G$ with $\mathfrak g_{\bar 0}:=Lie(G)$, existence of a universal complexification $\gamma:G \to G^\mathbb C$ is stated in \cite{Hoch}. In detail there is an ideal $\mathfrak p \subset \mathfrak g_{\bar 0}$ such that the Lie algebra $\mathfrak g^\mathbb C_{\bar 0}:=Lie(G^\mathbb C)$ is isomorphic to $(\mathfrak g_{\bar 0} \otimes  \mathbb C)/(\mathfrak p \otimes \mathbb C)$. The map $\gamma$ is given on Lie algebra level by the real embedding $\mathrm{emb}_{\bar 0}:\mathfrak g_{\bar 0} \to \mathfrak g_{\bar 0} \otimes  \mathbb C$ followed by projection. 
Approaching a  universal complexification for Lie supergroups, we find:
\begin{lemma}
 Let $(G,\mathfrak g_{\bar 0})$ be a real Harish-Chandra pair  and $\gamma:G\to G^\mathbb C$ its universal complexification. Set $\mathfrak g^\mathbb C_{\bar 0}:=Lie(G^\mathbb C)$. Let further $(G,\mathfrak g_{\bar 0}\oplus\mathfrak g_{\bar 1})$ be a real Harish-Chandra superpair. 
 Then  $\mathfrak g^\mathbb C:=\mathfrak g^\mathbb C_{\bar 0} \oplus (\mathfrak g_{\bar 1}\otimes \mathbb C)$ is a complex Lie superalgebra with respect to the inherited Lie superbracket. In particular  $(G^\mathbb C,\mathfrak g^\mathbb C)$ is a complex Harish-Chandra superpair.
\end{lemma}
\begin{proof}
 The adjoint representation of $\mathfrak g_{\bar 0}$ on $\mathfrak g_{\bar 1}$ integrates to an action $G\to GL_\mathbb C(\mathfrak g_{\bar 1}\otimes \mathbb C)$. The universal complexification yields $G^\mathbb C\to GL_\mathbb C(\mathfrak g_{\bar 1}\otimes \mathbb C)$. So on Lie superalgebra level $\mathfrak p\otimes \mathbb C$ acts trivially on $\mathfrak g_{\bar 1}\otimes \mathbb C$, i.e.~it is an ideal in $\mathfrak g\otimes\mathbb C$. This proves the lemma.
\end{proof}

\begin{theorem}
 Let $\mathcal G$ be a real Lie supergroup  associated to the Harish-Chandra superpair $(G,\mathfrak g)$ and let $\gamma:G\to G^\mathbb C$ be the universal complexification of $G$. Then the complex Harish-Chandra superpair $(G^\mathbb C, \mathfrak g^\mathbb C)$ together with the morphism $(\gamma,\Gamma_\ast)$, $\Gamma_\ast:=D_e\gamma \oplus \mathrm{emb}_{\bar 1} $ is associated to  a universal complexification $\mathcal G^\mathbb C$ of $\mathcal G$.
\end{theorem}
\begin{proof}
 Let $(H,\mathfrak h)$ be a complex Harish-Chandra superpair and $(\varphi,\Phi_\ast):(G,\mathfrak g)\to (H,\mathfrak h)$ be a morphism of real Harish-Chandra superpairs. Let $\varphi^\mathbb C:G^\mathbb C \to H$ be the underlying complexification and set $\Phi_\ast^\mathbb C:=D_e\varphi^\mathbb C \oplus \sigma$, where $\sigma:\mathfrak g_{\bar 1}\otimes \mathbb C \to \mathfrak h_{\bar 1}$ is the complex linear continuation of $\Phi_\ast|_{\mathfrak g_{\bar 1}}$. Then $(\varphi^\mathbb C, \Phi_\ast^\mathbb C)$ is unique with the required properties.
\end{proof}

\end{document}